\documentclass[11pt,a4paper]{article}
\usepackage[margin=25mm,top=23mm,bottom=23mm]{geometry}
\usepackage{amsmath,amssymb,amsthm}
\usepackage[hidelinks]{hyperref}
\hypersetup{pdftitle={Piecewise-polynomial invariants of rank-four matrix mutations},
 pdfsubject={A proof by returns of alternating mutations}}
\newtheoremstyle{result}{8pt}{4pt}{\normalfont}{}{\bfseries}{.}{.55em}{}
\theoremstyle{result}
\newtheorem{thm}{Theorem}[section]
\newtheorem{lem}[thm]{Lemma}
\newtheorem{defn}[thm]{Definition}
\newtheorem{cor}[thm]{Corollary}
\newcommand{\R}{\mathbb R}
\newcommand{\A}{\mathcal A}
\newcommand{\Int}{\operatorname{int}}
\newcommand{\Pf}{\operatorname{Pf}}
\newcommand{\pos}[1]{[#1]^+}
\newcommand{\T}{^{\mathsf T}}
\newcommand{\artlink}[2]{\href{#1}{#2}}
\numberwithin{equation}{section}
\begin{document}

\title{The description of the invariants of $4\times 4$ skew-symmetric matrix under cluster mutations.}
\author{G.~Chelnokov {\small\tt grishabenruven@yandex.ru }
 }

\date{}
\maketitle \maketitle

\begin{abstract}
We prove that any piecewise
polynomial invariant of a 4-quiver in each piece belongs to the ring, generated by Pfaffian and the determinant of symmetric companion. See precise formulation in Theorem 1.1.


\end{abstract}

\section{Statement of the result}\label{sec:main}
Fomin and Zelevinsky introduced clusters and cluster algebras, together with the matrix mutation that governs the change of an exchange matrix \cite{FZ}. An overview of the remarkable results that have followed since then is beyond the scope of this article. We study the mutation of a real skew-symmetric exchange matrix \(X=(x_{ij})\) of order four. Thus \(x_{ii}=0,\ x_{ji}=-x_{ij}\), and the coordinate ring and Pfaffian are
\[
 \A=\R[a,b,c,d,e,f],\quad
 (a,b,c,d,e,f)=(x_{12},x_{13},x_{14},x_{23},x_{24},x_{34}),\quad
 \pi=af-be+cd .
\]
The author previously proved (\cite[Theorem~1]{Ch}) that an invariant given by one polynomial on each orthant is a polynomial in \(\det X=\pi^2\), but the restriction of at most one piece at each orthant seems unnatural. Here we allow any finite partition into pieces \(S\) satisfying \(S\subseteq\overline{\Int S}\), and attach a polynomial \(P_S\) to each piece.

An open orthant determines an orientation of the complete graph: \(i\to j\) means \(x_{ij}>0\). A \emph{vortex} consists of a directed triangle and a fourth vertex that is a source or a sink. We call all other orientations \emph{admissible}. After refinement by orthants, let \(\Gamma(S)\) be the orientation of the interior of \(S\). Two pieces are joined in the \emph{transition graph} if a mutation takes an interior point of one to an interior point of the other (Definition~\ref{def:graph}).

\begin{thm}[Classification of polynomial branches]\label{th:main}
Let \(\mathcal C\) be a connected component of the transition graph of a finite-piece polynomial invariant.
\begin{enumerate}
\item If \(\mathcal C\) contains a piece with a vortex orientation, then \(P_S\in\R[\pi]\) for every \(S\in\mathcal C\).
\item If every orientation in \(\mathcal C\) is admissible, then \(P_S\in\R[\pi,K_{\Gamma(S)}]\) for every \(S\in\mathcal C\). Here
\[
 K_\Gamma=(\pi^2+16-\det C_\Gamma)/4,
\]
where \(C_\Gamma\) is the symmetric companion of Definition~\ref{def:comp}.
\end{enumerate}
One polynomial determines all branches on a component. Across each edge its Pfaffian argument changes sign and its \(K\)-argument, when present, is preserved. If the component has a closed walk of odd length, the polynomial is even in its Pfaffian argument.
\end{thm}

Fix a starting branch and two vertices. Two visits to the same target piece give a polynomial symmetry of the starting branch. Applying this to all six pairs yields six differential equations on that one branch. We determine their polynomial solutions by a rank calculation and an explicit polynomial section. Section~\ref{sec:finish} also identifies the system on the two admissible orientations with the Poisson system studied by Bondal.

\clearpage
\section{Open transitions and polynomial identities}\label{sec:graph}
\begin{defn}[Mutation and branch]\label{def:mut}
Write \(\pos{r}=\max(r,0)\). Mutation \(\mu_k\) negates all entries incident to \(k\). For \(u,v\ne k\) its formula is
\begin{equation}\label{eq:mut}
 (\mu_k X)_{uv}=x_{uv}+\pos{x_{uk}}\pos{x_{kv}}-\pos{-x_{uk}}\pos{-x_{kv}} .
\end{equation}
This is the real skew-symmetric form of matrix mutation in \cite[Definition~4.2]{FZ}. Suppose \(\R^6=\bigsqcup S\) is a finite partition into nonempty pieces with \(S\subseteq\overline{\Int S}\). A family of polynomials \(P_S\in\A\), called \emph{branches}, defines \(F(X)=P_S(X)\) for \(X\in S\). The function \(F\) is an \emph{invariant} if \(F(\mu_kX)=F(X)\) for all \(X\) and all \(k\).
\end{defn}
Each mutation is continuous and involutive: on the second application both incident entries in the increment change sign. Thus \(\mu_k\) is a homeomorphism. Its formula on each open orthant extends to a polynomial automorphism of \(\R^6\), since the incident entries can be recovered first and then the increments subtracted.

Every partition under consideration can be refined by orthants. For a piece \(S\), let \(V_\Gamma=\Int S\cap O_\Gamma\), where \(O_\Gamma\) runs over the open orthants. The union of the \(V_\Gamma\) is dense in \(\Int S\). There are finitely many of them, so every point of \(S\) belongs to the closure of a nonempty \(V_\Gamma\). Assign each point outside their union to one such \(V_\Gamma\). Each resulting piece contains \(V_\Gamma\), is contained in its closure and in \(S\), and has interior \(V_\Gamma\). It inherits the polynomial \(P_S\).

\begin{defn}[Transition graph]\label{def:graph}
The vertices are the refined pieces. An edge labelled \(k\) joins \(S\) to \(S'\) when \(\Int S\cap\mu_k^{-1}(\Int S')\ne\varnothing\). The graph is undirected because \(\mu_k\) is an involution. The definition allows loops and multiple edges.
\end{defn}
\begin{lem}[Transport along an edge or path]\label{le:transport}
If \(m\) is the polynomial formula on \(S\) for an edge \(S\to S'\), then \(P_S=P_{S'}\circ m\) in \(\A\). For a path \(\gamma\) with successive formulas \(m_1,\ldots,m_r\), the same identity holds with \(m_\gamma=m_r\circ\cdots\circ m_1\).
\end{lem}
\begin{proof}
Invariance gives the equality on the nonempty open set in Definition~\ref{def:graph}. A polynomial vanishing there vanishes identically, by the one-variable identity theorem applied successively in each coordinate on a small open box. Compose the edge identities to obtain the path identity. The image of a single piece can meet the interiors of several different pieces.
\end{proof}
\begin{lem}[Sign of the Pfaffian]\label{le:pf}
For every mutation, \(\pi(\mu_kX)=-\pi(X)\).
\end{lem}
\begin{proof}
For the alternating form with matrix \(X\), set \(e'_k=-e_k\) and \(e'_i=e_i+\pos{x_{ik}}e_k\) for \(i\ne k\). Its matrix in the new basis is \eqref{eq:mut} since the increment equals \(\pos{x_{ik}}x_{kj}+x_{ik}\pos{x_{jk}}\). The basis change has determinant \(-1\). For \(\omega=\sum_{i<j}x_{ij}e_i^*\wedge e_j^*\), the value of \(\omega\wedge\omega/2\) on the ordered basis is \(af-be+cd\). It is multiplied by the determinant of the basis change.
\end{proof}
An odd path may return to a piece containing both signs of \(\pi\). If one further splits the interiors along \(\pi=0\), the two signs color the refined transition graph, which is then bipartite by Lemma~\ref{le:pf}. Boundary points can be assigned by closure as in the orthant refinement.

\clearpage
\section{The second generator on admissible orientations}\label{sec:comp}
\begin{defn}[Symmetric companion]\label{def:comp}
For \(X\) in an open orthant a companion is a symmetric matrix \(C=(c_{ij})\) satisfying \(c_{ii}=2,\ |c_{ij}|=|x_{ij}|\) for \(i\ne j\). We require \(c_{ij}c_{jk}c_{ki}>0\) on a directed triangle and \(c_{ij}c_{jk}c_{ki}<0\) on a transitive triangle.
\end{defn}
These signs are Seven's admissibility condition for a quasi-Cartan companion \cite[Definition~2.10]{Seven}. In the complete graph the chordless cycles are exactly the triangles. We give the argument for real entries.
\begin{lem}[Existence of the companion]\label{le:comp}
A companion exists exactly on the admissible orientations. It is unique up to simultaneous sign changes of the corresponding rows and columns. In particular, \(K_\Gamma\) is a well-defined element of \(\A\) for every admissible orthant.
\end{lem}
\begin{proof}
An acyclic orientation has no directed triangles. If there is a directed triangle, a fourth vertex that is a source or a sink gives no further directed triangle. Otherwise it points to either one or two vertices of the triangle. Inspection of the three cyclically ordered edges gives exactly one additional directed triangle. The latter orientation is strongly connected, and the cyclic symmetry of the triangle and reversal of all arrows give all such orientations up to relabelling and reversal.

The product of the four prescribed signs must be positive, since every edge occurs twice. This condition is equivalent to having an even number of directed triangles. Conversely, choose the signs on the three edges incident to vertex 1. Their three triangles determine the remaining edge signs, and the fourth condition follows by parity. Sign changes of basis vectors normalize the first three signs, proving uniqueness. They leave the determinant unchanged. On a fixed orthant \(c_{ij}\) is a fixed signed copy of \(x_{ij}\) for \(i<j\), proving polynomiality.
\end{proof}
\begin{lem}[Transport of the companion]\label{le:K}
If an open mutation transition leads from an admissible orientation \(\Gamma\) to another admissible orientation \(\Gamma'\), then \(K_{\Gamma'}\circ m=K_\Gamma\).
\end{lem}
\begin{proof}
Use the companion basis change of Barot, Geiss and Zelevinsky \cite[Proposition~3.2]{BGZ}. Put \(e'_k=-e_k\), put \(e'_i=e_i-c_{ik}e_k\) for \(i\to k\), and keep \(e'_j=e_j\) for \(k\to j\). Its determinant is \(-1\). The diagonal entries remain 2. The coefficient between two incoming or two outgoing vertices is unchanged. A coefficient incident to \(k\) stays the same on the incoming side and changes sign on the outgoing side.

For a path \(i\to k\to j\), write \(t=c_{ik}c_{kj}\). Its new outer coefficient is \(c_{ij}-t\) and the new triangle product is \(t^2-tc_{ij}\). In a transitive triangle \(tc_{ij}<0\): the magnitudes add and the new triangle is directed. In a directed triangle \(tc_{ij}>0\): the new magnitude is the absolute difference and the new triangle is directed exactly when \(|t|>|c_{ij}|\). Equality would put the image on a coordinate hyperplane. These are the alternatives in \eqref{eq:mut}. For two incoming or two outgoing vertices the triangle remains transitive and its product remains negative. This checks all three triangles through \(k\). Admissibility at the target forces the fourth sign by the parity argument in Lemma~\ref{le:comp}. The new companion has the same determinant. Lemmas~\ref{le:transport} and~\ref{le:pf}, or polynomial continuation from the open transition, now give the claim.
\end{proof}

\clearpage
\section{Alternating mutations at two vertices}\label{sec:pair}
Fix a pair and order its vertices so that \(h=x_{ij}>0\) on the starting piece. Write \(k<\ell\) for the other two vertices. Set
\begin{equation}\label{eq:chart}
 p_s=x_{is},\quad q_s=x_{js}\quad(s=k,\ell),\qquad
 w=x_{k\ell},\qquad \Pi=hw-p_kq_\ell+p_\ell q_k .
\end{equation}
The sign of the permutation \((i,j,k,\ell)\) determines whether \(\Pi=\pi\) or \(-\pi\). The six quantities \(h,\Pi,p_k,q_k,p_\ell,q_\ell\) are independent coordinates after \(h\) is inverted. Equation~\eqref{eq:chart} recovers \(w\) by division by \(h\).

\begin{lem}[One step on the two columns]\label{le:step}
Mutate successively at \(i,j,i,j,\ldots\). After each individual step exchange the names of the first and second vertices of the pair. Then the positive entry \(h\) is unchanged. On each column, in coordinates \((u,v)=(q_s,-p_s)\), the step is
\begin{equation}\label{eq:T}
 T_h(u,v)=(v,-u+h\pos{v}).
\end{equation}
After an even number of steps the original labels are restored and both \(h\) and \(\Pi\) have their initial values.
\end{lem}
\begin{proof}
Before exchanging the names, mutation at the first vertex sends \(p_s\) to \(-p_s\) and \(q_s\) to \(q_s-h\pos{-p_s}\) by \eqref{eq:mut}. Exchanging the names gives \eqref{eq:T}. Both the mutation and the exchange negate the entry of the chosen pair, so the next step uses the same positive \(h\). The assertion about \(\Pi\) follows from Lemma~\ref{le:pf}.
\end{proof}
Exchanging the two plane coordinates identifies \(T_h\) with the map \(T_{h,0}\) in the family studied by Lagarias and Rains \cite{LR}. Here the same map acts on two columns. Their relative timing under iteration is the point of Lemma~\ref{le:sync}.

\begin{lem}[Powers of one matrix]\label{le:power}
Write
\[
 J=\begin{pmatrix}0&1\\-1&0\end{pmatrix},
 \qquad H=\begin{pmatrix}0&1\\-1&h\end{pmatrix}.
\]
For a column \(z=(u,v)\T\), put \(\widehat z=J^{-1}z\) in the third quadrant \(u<0,v<0\), and \(\widehat z=z\) elsewhere. Along a trajectory avoiding the axes, after \(m\) steps one has
\begin{equation}\label{eq:power}
 \widehat z_m=\sigma(m)H^{r(m)}\widehat z_0,\qquad
 r(m)=\#\{0\le n<m:v_n>0\},\qquad \sigma(m)\in\{1,-1\}.
\end{equation}
\end{lem}
\begin{proof}
If \(v>0\), \eqref{eq:T} acts by \(H\), and the image has positive first coordinate. If \(u>0,v<0\), it acts by \(J\) and the image is in the third quadrant, so the normalized column is unchanged. If \(u<0,v<0\), the image \(Jz\) is in the second quadrant. The step in normalized coordinates is then \(J^2=-I\). Multiplying these three possible step matrices proves \eqref{eq:power}.
\end{proof}
Thus the two planar columns determine the entire matrix at every even time: use the constant \(\Pi\) in \eqref{eq:chart} to recover the sixth entry. The normalization is a linear change on each fixed quadrant; it will always be used with the quadrant specified.

\clearpage
\section{Synchronized powers and two visits to one piece}\label{sec:sync}
\begin{lem}[Synchronization of the exponents]\label{le:sync}
Suppose two columns have the same \(h>0\), and both trajectories avoid the axes. Their exponents in \eqref{eq:power} satisfy
\[
 r_s(m)\ge\lfloor m/3\rfloor\quad(s=k,\ell),\qquad
 |r_k(m)-r_\ell(m)|\le4 \qquad(m\ge0).
\]
\end{lem}
\begin{proof}
A column in the fourth quadrant next enters the third, and one in the third next enters the second. Hence \(v\) cannot be negative three times in succession. This gives the lower bound. Let \(b(m)\) and \(c(m)\) count visits at times \(0,\ldots,m-1\) to the third and fourth quadrants, respectively. Every visit to the fourth is followed by one to the third; every visit to the third has come from the fourth, except possibly at time zero. Thus \(|b(m)-c(m)|\le1\).

Compare the two counts \(b(m)\) by using directions. The positively homogeneous map \(T_h\) induces an orientation-preserving circle homeomorphism: its inverse is \((U,V)\mapsto(h\pos U-V,U)\), and both linear parts have determinant 1. Represent the initial angle by \(\alpha\in[\pi,3\pi)\); in the angle formulas \(\pi\) denotes the circle constant. Choose the increasing lift \(G\) of the angle map with \(G(\pi)=\pi/2\) and \(G(\alpha+2\pi)=G(\alpha)+2\pi\). Then \(G(3\pi/2)=\pi\) and \(G(3\pi)=5\pi/2\). On reducing the image angle to the chosen interval, one adds \(2\pi\) exactly for a step starting in the third quadrant. Iteration gives
\[
 b(m)=-\left\lfloor\frac{G^m(\alpha)-\pi}{2\pi}\right\rfloor.
\]
Two initial angles in this interval remain less than \(2\pi\) apart under the increasing lift \(G^m\), so \(|b_k-b_\ell|\le1\). Finally \(r=m-b-c=m-2b+(b-c)\), giving the bound 4.
\end{proof}

\begin{lem}[Return identity on the starting branch]\label{le:return}
For every branch \(P_S\) and every pair \(i,j\), there exists an integer \(d>0\) such that, in the starting coordinates \((h,\Pi,\widehat z_k,\widehat z_\ell)\), the branch is invariant under simultaneous multiplication of both columns by \(H^d\), with \(h,\Pi\) fixed.
\end{lem}
\begin{proof}
Choose a point of \(\Int S\) all of whose iterates belong to piece interiors and avoid the coordinate hyperplanes. Such points are dense: the union of interiors is open and dense by the hypothesis on the pieces, and remains so after removing the hyperplanes. Its inverse images under finite mutation words are open and dense. Their countable intersection is dense, as is seen by choosing nested closed balls of radii tending to zero inside any given open ball and successively inside these open sets.

Among the even times choose an infinite subsequence with the same target piece, both target quadrants, both signs \(\sigma_s\), and the same difference \(r_k-r_\ell\). There are finitely many choices by Lemmas~\ref{le:power} and~\ref{le:sync}. Take two of these times at which \(r_k\) increases by \(d>0\). The increase in \(r_\ell\) is also \(d\). Each finite itinerary has a fixed formula on an open neighborhood of the chosen point. In the chosen source and target coordinates the two formulas \(A_1,A_2\) satisfy \(A_2=A_1\circ B\), where \(B\) applies \(H^d\) to both columns and fixes \(h,\Pi\). For the common target branch \(Q\), Lemma~\ref{le:transport} yields \(P_S=Q\circ A_1=Q\circ A_2=P_S\circ B\). These are identities in the ring with \(h^{-1}\) adjoined. The target piece can differ from \(S\).
\end{proof}

\clearpage
\section{Six equations on a single polynomial}\label{sec:D}
\begin{lem}[From a return to a derivation]\label{le:weight}
Suppose \(P\) is a polynomial in the two normalized columns over \(\R(h,\Pi)\) and is invariant under simultaneous multiplication by \(H^d\) for some \(d>0\). Then \(DP=0\), where \(D\) fixes \(h,\Pi\) and acts on each column by \(2H-hI\).
\end{lem}
\begin{proof}
A derivation is a linear operator satisfying \(D(fg)=fDg+gDf\). Its values on the coordinates determine its action. Adjoin a root \(\lambda\) of \(\lambda^2-h\lambda+1=0\) to the coefficient field; the other root is \(\lambda^{-1}\). They are distinct, and \(\lambda\) is not a root of unity since \(h=\lambda+\lambda^{-1}\) is an indeterminate.

For each normalized column \((U,V)\) take the independent forms \(\xi=U-\lambda V\) and \(\eta=U-\lambda^{-1}V\). Under \(H\) these acquire the factors \(\lambda,\lambda^{-1}\), respectively. On a monomial in the four forms the factor under \(H^d\) is \(\lambda^{dm}\), where \(m\) is its total \(\xi\)-degree minus its total \(\eta\)-degree. Invariance leaves only monomials with \(m=0\). The derivation induced by \(2H-hI\) multiplies the two forms by opposite scalars \(\pm(\lambda-\lambda^{-1})\). By the product rule it multiplies the monomial by \(m(\lambda-\lambda^{-1})\), which proves the claim.
\end{proof}
\begin{lem}[The six derivations]\label{le:D}
In a fixed orthant order each of the six pairs by its arrow \(i\to j\). In the coordinates \eqref{eq:chart} put \(\varepsilon_s=-1\) if \(p_s>0,\ q_s<0\), and \(\varepsilon_s=1\) otherwise. Every branch \(P_S\) in this orthant satisfies \(D_{ij}P_S=0\), where
\begin{equation}\label{eq:D}
\begin{aligned}
 D_{ij}h&=0,\\
 D_{ij}p_s&=\varepsilon_s h p_s+2q_s,&
 D_{ij}q_s&=-2p_s-\varepsilon_s hq_s\quad(s=k,\ell),\\
 D_{ij}w&=(\varepsilon_k-\varepsilon_\ell)(p_kq_\ell+q_kp_\ell).
\end{aligned}
\end{equation}
\end{lem}
\begin{proof}
Apply Lemmas~\ref{le:return} and~\ref{le:weight} separately to each pair on the chosen starting piece. In the original column \((q_s,-p_s)\T\) the matrix is \(2H-hI\) outside the third quadrant and \(J(2H-hI)J^{-1}\) inside it. This gives the two column formulas. Differentiating \(\Pi\) in \eqref{eq:chart} and using \(D_{ij}\Pi=0\) gives the last formula after cancelling \(h\). All six coordinate values are polynomial, so the identity obtained in the field holds in \(\A\).
\end{proof}
\begin{cor}[The ring for one pair]\label{cor:one}
Put \(F_0=\R(h,\Pi)\), and write \((U_s,V_s)\T=\widehat z_s\). In the polynomial ring on the two columns over \(F_0\), the kernel of \(D_{ij}\) is generated by
\[
 Q_{kk},\quad Q_{k\ell},\quad Q_{\ell\ell},\quad
 \Delta=U_kV_\ell-U_\ell V_k,\qquad
 Q_{st}=U_sU_t-\tfrac h2(U_sV_t+V_sU_t)+V_sV_t .
\]
\end{cor}
\begin{proof}
This is the two-vector case of the first fundamental theorem for the special orthogonal group of a binary quadratic form \cite[Section~10.2(b)]{KP}. For a direct proof over \(F_0\), use the weights from Lemma~\ref{le:weight}. After adjoining \(\lambda\), a weight-zero monomial factors into products \(\xi_s\eta_t\). The diagonal products are \(Q_{ss}\); the mixed products have sum \(2Q_{k\ell}\) and difference \((\lambda-\lambda^{-1})\Delta\). All four stated generators have weight zero. Averaging the coefficients of a representation with their conjugates under \(\lambda\leftrightarrow\lambda^{-1}\) gives a representation over \(F_0\).
\end{proof}

\clearpage
\section{Rank and algebraic independence}\label{sec:rank}
The rank of a family of derivations is the rank, over \(E=\R(a,b,c,d,e,f)\), of the matrix of their values on \(a,b,c,d,e,f\). It is the largest order of a minor that is a nonzero polynomial.
\begin{lem}[Rank on the three orientations]\label{le:rank}
System~\eqref{eq:D} has rank at least four on every orthant and at least five on a vortex.
\end{lem}
\begin{proof}
The linear terms in \eqref{eq:D} are \(Dp_s=2q_s,\ Dq_s=-2p_s,\ Dh=Dw=0\). Reversing the order of the endpoints of a pair negates the linear part of that row. Use the rows for \(13,14,23,24\) and the columns for \(b,c,d,e\). At \(a=1,\ f=2,\ b=c=d=e=0\) the matrix of linear parts, up to those row signs, is
\[
 \begin{pmatrix}0&4&2&0\\-4&0&0&2\\-2&0&0&4\\0&-2&-4&0\end{pmatrix},
 \qquad \det=144 .
\]
The degree-four part of the corresponding minor is the determinant of the linear parts, so the minor is a nonzero polynomial. The evaluation point can lie outside the orthant for this polynomial test.

All vortices are related by relabelling and reversal of every arrow. Relabelling permutes indices in \eqref{eq:D}. Under \(X\mapsto-X\) the order of each pair is reversed, \((p_s,q_s)\) becomes \((-q_s,-p_s)\), and \(\varepsilon_s\) stays the same. Substitution in \eqref{eq:D} negates the derivation after the coordinate change and preserves the rank. Consider the vortex \(1\to2\to3\to1\) with sink 4. At \((a,b,c,d,e,f)=(1,-2,3,4,5,6)\), the rows \(12,31,14,23,24\) and columns \(a,b,c,d,e\) in \eqref{eq:D} give
\[
 \begin{pmatrix}
 0&6&13&0&-11\\6&0&-18&6&0\\-13&-18&0&-8&17\\
 0&6&0&0&32\\11&-12&-17&-32&0
 \end{pmatrix},\qquad \det=48384 .
\]
This proves the second assertion.
\end{proof}
\begin{lem}[Rank bounds the number of independent functions]\label{le:diff}
If \(r\) derivations on \(E=\R(x_1,\ldots,x_n)\) are linearly independent over \(E\), their common kernel contains at most \(n-r\) algebraically independent functions over \(\R\).
\end{lem}
\begin{proof}
The gradients of elements of the common kernel annihilate the \(r\) independent rows, so they lie in a subspace of dimension \(n-r\). Algebraically independent functions \(g_1,\ldots,g_m\) have independent gradients over \(E\). To see this, extend the \(g_i\) to a transcendence basis by adjoining each \(x_j\) that is still transcendental over the field constructed so far. On the rational field in this basis the partial derivatives with respect to \(g_i\) exist. They extend to \(E\): for an algebraic element \(v\) with minimal polynomial \(b(T)\), set \(Dv=-(Db)(v)/b'(v)\), differentiating only the coefficients in \(Db\). The denominator is nonzero in characteristic zero. This rule makes the relation \(b(v)=0\) stable and so defines a derivation on the extension. Pairing a linear dependence of the gradients with each of these extended derivations makes its coefficients vanish one at a time.
\end{proof}

\clearpage
\section{Polynomial sections}\label{sec:sections}
\begin{lem}[From algebraic dependence to a polynomial formula]\label{le:section}
Let \(g=(g_1,\ldots,g_m)\colon\R^n\to\R^m\) be polynomial and suppose it admits a polynomial section, that is, a polynomial map \(L\colon\R^m\to\R^n\) such that \(g(L(s))=s\). Then \(g_1,\ldots,g_m\) are algebraically independent. Every polynomial \(P\) algebraic over \(\R(g_1,\ldots,g_m)\) satisfies
\begin{equation}\label{eq:section}
 P=(P\circ L)(g_1,\ldots,g_m)\in\R[g_1,\ldots,g_m].
\end{equation}
\end{lem}
\begin{proof}
Substituting \(x=L(s)\) in any relation among the \(g_i\) gives the same relation among the independent \(s_i\). Identify \(\R(g)\) with \(\R(s)\) and let \(b(s,T)\) be the monic minimal polynomial of \(P\) over this field. Clear the denominators in \(b(g,P)=0\) before substituting \(x=L(s)\). The common denominator remains a nonzero polynomial in \(s\). Hence \(b(s,P(L(s)))=0\) over \(\R(s)\). An irreducible polynomial with a root in its ground field is linear. Its root is \(P(L(s))\), which proves \eqref{eq:section}.
\end{proof}
\begin{lem}[Sections for the generators]\label{le:explicit}
For every admissible orthant the polynomial map \((\pi,K_\Gamma)\) admits a polynomial section. The Pfaffian alone also admits a polynomial section.
\end{lem}
\begin{proof}
For \(\pi\) take \(a=1,\ f=s,\ b=c=d=e=0\). For the pair, take two representatives as in Lemma~\ref{le:comp}: the positive acyclic orthant, and the strongly connected orthant with signs \((-,+,+,-,-,+)\) in the order \(a,b,c,d,e,f\). With \(\Sigma=a^2+b^2+c^2+d^2+e^2+f^2\), their formulas are
\begin{equation}\label{eq:K}
\begin{aligned}
 K_{\rm ac}&=\Sigma+abd+ace+bcf+def+acdf,\\
 K_{\rm sc}&=\Sigma-abd-ace+bcf+def-abef.
\end{aligned}
\end{equation}
The strictly upper entries of the two companions may be chosen as \((-a,-b,-c,-d,-e,-f)\) and \((a,-b,-c,-d,-e,-f)\), respectively. Checking the four signs in Definition~\ref{def:comp} and expanding \((\pi^2+16-\det C)/4\) gives \eqref{eq:K}.

Let \(\epsilon=1\) in the first case and \(-1\) in the second. Set \(c=f=0,\ b=5/2,\ e=-2s/5\). Then \(\pi=s\) and
\[
 K=b^2+e^2+a^2+d^2+\tfrac52\epsilon ad
   =b^2+e^2+(a+2\epsilon d)(a+\epsilon d/2).
\]
Put \(V=t-b^2-e^2\), \(a=(4V-1)/3\), and \(d=2\epsilon(1-V)/3\). The two factors are then \(1,V\), giving \(K=t\). A section is required for the polynomial formulas on the whole affine space, so its values may lie outside the defining orthant. Transport the section to any other representative by relabelling and global arrow reversal. These operations preserve the companion construction. An odd relabelling also changes the sign of the section parameter \(s\).
\end{proof}
The first expression in \eqref{eq:K} coincides with the Markov invariant \(M_Q\) of Fomin and Neville for the indicated vertex order \cite[Example~9.3, (9.6)]{FN}; their \((x,y,u,z,v,w)\) equal our \((a,d,f,b,e,c)\). It also equals \(C_1+4\) for the Stokes-matrix Casimir in Ugaglia \cite[Example~2]{Ug}, upon setting her \((p,q,r,x,y,z)=(-a,-b,-c,-d,-e,-f)\).

\clearpage
\section{Intersection of the kernels and transport}\label{sec:finish}
\begin{proof}[Proof of Theorem~\ref{th:main}]
Start with a vortex piece \(S\). Let \(E_0\) be the common kernel of its six derivations in \(E=\R(a,b,c,d,e,f)\). It is a field by the quotient rule. It contains \(\pi\) by \eqref{eq:chart} and \eqref{eq:D}, and \(P_S\) by Lemma~\ref{le:D}. Lemmas~\ref{le:rank} and~\ref{le:diff} show that \(E_0\) has transcendence degree at most one. The nonconstant \(\pi\) is transcendental over \(\R\), so \(P_S\) is algebraic over \(\R(\pi)\). Lemmas~\ref{le:explicit} and~\ref{le:section} give \(P_S=R(\pi)\) with \(R\in\R[t]\). For any neighbor \(S'\), apply Lemmas~\ref{le:transport} and~\ref{le:pf} to the reverse edge to obtain \(P_{S'}=R(-\pi)\). Following paths covers the whole component.

Now suppose the component is entirely admissible and fix a piece \(S\) in it. For each pair, \(D_{ij}K_{\Gamma(S)}=0\). Indeed, by Lemma~\ref{le:K} the family \(K_{\Gamma(S')}\) obeys the polynomial transport identity on every edge in the component. For the point chosen in Lemma~\ref{le:return} all intermediate steps are open transitions, so all its pieces belong to this component. The proofs of Lemmas~\ref{le:return}, \ref{le:weight} and~\ref{le:D} therefore apply to this family. The common kernel \(E_0\) for \(S\) contains \(\pi,K_{\Gamma(S)},P_S\). Its transcendence degree is at most two by Lemmas~\ref{le:rank} and~\ref{le:diff}. Its first two displayed elements are algebraically independent by Lemmas~\ref{le:section} and~\ref{le:explicit}. The third is therefore algebraic over \(\R(\pi,K_{\Gamma(S)})\). The same section gives
\[
 P_S=R_S(\pi,K_{\Gamma(S)}),\qquad R_S\in\R[t,u].
\]

Across an edge \(S\to S'\), Lemmas~\ref{le:transport}, \ref{le:pf} and~\ref{le:K} and algebraic independence give \(R_{S'}(t,u)=R_S(-t,u)\). The same formula in one variable holds on a component containing a vortex. Along a path of length \(r\) the sign is \((-1)^r\). In particular, an odd closed walk forces evenness in the first argument.
\end{proof}

\begin{cor}[Transported rings]\label{cor:paths}
Fix a piece \(S_0\). Suppose \(P_{S_i}\) lies in a known ring \(A_i\subseteq\A\) for \(1\le i\le r\). For arbitrary paths \(\gamma_i\colon S_0\to S_i\),
\[
 P_{S_0}\in\bigcap_{i=1}^r m_{\gamma_i}^*A_i,\qquad
 m_{\gamma_i}^*Q=Q\circ m_{\gamma_i}.
\]
\end{cor}
\begin{proof}
Apply Lemma~\ref{le:transport} along each path.
\end{proof}
In Lemma~\ref{le:return} the two paths for a given pair have the same source and the same target. Comparing them gives the equation on \(P_{S_0}\) directly. The argument above computes the common polynomial kernel as \(\R[\pi]\) on a vortex and as \(\R[\pi,K_\Gamma]\) on each piece of an entirely admissible component. On other pieces of a component containing a vortex, transport from the vortex enforces the stronger restriction to \(\R[\pi]\).

\clearpage
\subsection*{The admissible kernels and Bondal's Poisson center}
There is a second way to identify the common polynomial kernel on an admissible orthant. Bondal gives the coordinate Poisson bracket on the space of upper unitriangular matrices in \cite[(5.3)]{Bondal}, and describes its full polynomial center in \cite[Theorem~5.1]{Bondal}. Ugaglia also obtained this bracket from Stokes matrices \cite{Ug}. For a Poisson bracket, the \emph{center} consists of all polynomials whose bracket with every coordinate vanishes. We give the coordinate identification to specify where this theorem applies.

\begin{lem}[The system in Bondal's coordinates]\label{le:bondal}
Order the vertices so that every arrow goes from a smaller to a larger index in the acyclic case; in the strongly connected case take the order with the sole exception \(4\to1\). Put \(U_{ii}=1,\ U_{ij}=-x_{ij}\) for \(i<j\), and \(U_{ij}=0\) for \(i>j\). With Bondal's normalization of the bracket, on each directed edge \(i\to j\)
\[
 D_{ij}=-2\{x_{ij},\,\cdot\,\},\qquad x_{ji}=-x_{ij}.
\]
Consequently, the common polynomial kernel on either admissible orientation is \(\R[\pi,K_\Gamma]\).
\end{lem}
\begin{proof}
For \(i<j<k\), formula~(5.3) of \cite{Bondal}, pulled back by \(U_{ij}=-x_{ij}\), reads
\begin{align*}
 \{x_{ij},x_{ik}\}&=-x_{jk}-\tfrac12x_{ij}x_{ik},&
 \{x_{ij},x_{jk}\}&= x_{ik}+\tfrac12x_{ij}x_{jk},\\
 \{x_{ik},x_{jk}\}&=-x_{ij}-\tfrac12x_{ik}x_{jk}.
\end{align*}
For \(i<j<k<\ell\), disjoint noncrossing pairs have bracket zero and
\(\{x_{ik},x_{j\ell}\}=x_{ij}x_{k\ell}-x_{jk}x_{i\ell}\).
In the acyclic order, \(\varepsilon_s=-1\) in \eqref{eq:D} exactly when \(i<s<j\). The two column formulas now give the three shared-index brackets above after multiplication by \(-2\); the formula for \(w\) gives the disjoint pairs. The strongly connected representative used in Lemma~\ref{le:explicit} takes the required order upon listing its vertices as \(1,3,4,2\). In the new numbering, reversing the extreme edge \(14\) leaves every \(\varepsilon_s\) for the five forward pairs unchanged: that edge occurs only in an exterior column, and the new signs in that column never have \(p_s>0,\ q_s<0\). For the remaining pair \(4\to1\) within each column \(p_s,q_s\) exchange roles. Both \(\varepsilon_s\) change from \(-1\) to \(+1\), and \(h\) is negated. Equation~\eqref{eq:D} gives \(D_{41}=-D^{\rm ac}_{14}\). This verifies the asserted identity in both cases.

Bondal's generators in even order are the coefficients of \(\chi(t)=\det(U+tU\T)\) and \(\Pf(U-U\T)\) \cite[Section~2.2, Theorem~5.1]{Bondal}. The theorem over \(\mathbb C\) gives the real statement by conjugating and averaging coefficients of a polynomial expression in the real generators. Here \(U-U\T=-X\), so \(\Pf(U-U\T)=\pi\). The symmetric matrix \(U+U\T\) is a companion in either chosen order: in the second order exactly the two triangles containing the backward edge are directed. Transposition gives the reciprocal identity \(\chi(t)=t^4\chi(1/t)\); write
\(\chi(t)=1+c_1t+c_2t^2+c_1t^3+t^4\). Hence
\[
 \chi(1)=\det C_\Gamma,\quad \chi(-1)=\pi^2,\qquad
 c_1=4-K_\Gamma,\quad c_2=\pi^2+2c_1-2 .
\]
This identifies the center with the stated polynomial ring.
\end{proof}
Bondal's theorem thus supplies the admissible kernel once Lemma~\ref{le:D} has placed a branch in it. The implication from finite-piece mutation invariance to the six equations uses Lemma~\ref{le:return}; the rank-five vortex calculation is in Lemma~\ref{le:rank}.

\section{Acknowledgements}
The role of the AI GPT-6 Astra in the creation of this paper is so significant that only the arxiv.org\` policy regarding the inclusion of AI as a 
co-author prevents as from listing GPT-6 Astra as one of authors.


\clearpage
\begin{thebibliography}{99}
\bibitem{BGZ}
M. Barot, C. Geiss and A. Zelevinsky,
\emph{Cluster algebras of finite type and positive symmetrizable matrices},
J. London Math. Soc. (2) \textbf{73} (2006), 545--564.
\artlink{https://arxiv.org/abs/math/0411341}{arXiv:math/0411341}.

\bibitem{Bondal}
A. I. Bondal, \emph{A symplectic groupoid of triangular bilinear forms and the braid group},
Izv. Math. \textbf{68} (2004), 659--708.
\artlink{https://www.mathnet.ru/eng/im495}{Math-Net.Ru: im495}.

\bibitem{Ch}
G. Chelnokov, \emph{On the invariants of \(4\times 4\) skew-symmetric matrix under cluster mutations} (2023).
\artlink{https://arxiv.org/abs/2311.01784}{arXiv:2311.01784}.

\bibitem{FN}
S. Fomin and S. Neville, \emph{Cyclically ordered quivers},
Adv. Math. \textbf{500} (2026), article 111059.
\artlink{https://arxiv.org/abs/2406.03604}{arXiv:2406.03604}.

\bibitem{FZ}
S. Fomin and A. Zelevinsky, \emph{Cluster algebras I: Foundations},
J. Amer. Math. Soc. \textbf{15} (2002), 497--529.
\artlink{https://arxiv.org/abs/math/0104151}{arXiv:math/0104151}.

\bibitem{KP}
H. Kraft and C. Procesi, \emph{Classical Invariant Theory: A Primer},
preliminary version (July 1996), Section~10.2(b).
\artlink{https://dmi.unibas.ch/fileadmin/user_upload/dmi/Personen/Kraft_Hanspeter/Classical_Invariant_Theory.pdf}{University of Basel edition}.

\bibitem{LR}
J. C. Lagarias and E. M. Rains, \emph{Dynamics of a family of piecewise-linear area-preserving plane maps I. Rational rotation numbers},
J. Difference Equ. Appl. \textbf{11} (2005), 1089--1108.
\artlink{https://arxiv.org/abs/math/0301294}{arXiv:math/0301294}.

\bibitem{Seven}
A. I. Seven, \emph{Cluster algebras and semipositive symmetrizable matrices}, Trans. Amer. Math. Soc. \textbf{363} (2011), 2733--2762.
\artlink{https://arxiv.org/abs/0804.1456}{arXiv:0804.1456}.

\bibitem{Ug}
M. Ugaglia, \emph{On a Poisson structure on the space of Stokes matrices},
Internat. Math. Res. Notices \textbf{1999}, no.~9, 473--493.
\artlink{https://arxiv.org/abs/math/9902045}{arXiv:math/9902045}.
\end{thebibliography}
\end{document}